\documentclass[leqno,12pt]{amsart}
\usepackage{amssymb} 
\theoremstyle{plain} 
\newtheorem{theorem}{\indent\sc Theorem}[section] 

\newtheorem{proposition}[theorem]{\indent\sc Proposition}

\theoremstyle{definition} 

\begin{document}

\title[Coupled Painlev\'e III system]{Coupled Painlev\'e III system with affine Weyl group symmetry of type $D_6^{(1)}$ \\}
\author{Yusuke Sasano }

\renewcommand{\thefootnote}{\fnsymbol{footnote}}
\footnote[0]{2000\textit{ Mathematics Subjet Classification}.
34M55; 34M45; 58F05; 32S65.}

\maketitle

\begin{abstract}
We find and study a six-parameter family of coupled Painlev\'e III systems in dimension six with affine Weyl group symmetry of type $D_6^{(1)}$. We also find and study its degenerate systems with affine Weyl group symmetry of types $B_5^{(1)}$ and $D_5^{(2)}$.
\end{abstract}

\section{Introduction}

In \cite{Sasa4,Sasa3}, we presented a 4-parameter family of 2-coupled Painlev\'e III systems in dimension four with affine Weyl group symmetry of type $D_4^{(1)}$. We will make non-linear ordinary differential systems with affine Weyl group symmetry of type $D_{2n+2}^{(1)} \ (n \geq 2)$. 

In \cite{Sasa1,Sasa5}, we succeeded to make $(2n+2)$-parameter family of n-coupled Painlev\'e VI systems in dimension 2n with affine Weyl group symmetry of type $D_{2n+2}^{(1)} \ (n \geq 1)$ by connecting the invariant divisors $p_i,q_i-q_{i+1},p_{i+1}$ for the canonical variables $(q_i,p_i) \ (i=1,2,\ldots,n)$. These systems are polynomial Hamiltonian systems with coupled Painlev\'e VI Hamiltonians given by
\begin{equation}
\frac{dq_i}{dt}=\frac{\partial H}{\partial p_i}, \quad \frac{dp_i}{dt}=-\frac{\partial H}{\partial q_i} \quad (i=1,2,\ldots,n)
\end{equation}
with the polynomial Hamiltonian
\begin{align}
\begin{split}
H &=\sum_{i=1}^n H_{VI}(q_i,p_i,t;\alpha_0^{(i)},\alpha_1^{(i)},\ldots,\alpha_4^{(i)})+\sum_{1 \leqslant l<m \leqslant n}\frac{2(q_l-t)p_lq_m((q_m-1)p_m+\alpha_2^{(m)})}{t(t-1)},
\end{split}
\end{align}
where the symbol $H_{VI}(x,y,t;\alpha_0,\alpha_1,\alpha_2,\alpha_3,\alpha_4)$ is given by
\begin{align}
\begin{split}
&H_{VI}(x,y,t;\alpha_0,\alpha_1,\alpha_2,\alpha_3,\alpha_4)\\
&=\frac{1}{t(t-1)}[y^2(x-t)(x-1)x-\{(\alpha_0-1)(x-1)x+\alpha_3(x-t)x\\
&+\alpha_4(x-t)(x-1)\}y+\alpha_2(\alpha_1+\alpha_2)x]  \quad (\alpha_0+\alpha_1+2\alpha_2+\alpha_3+\alpha_4=1). 
\end{split}
\end{align}
However, in this case the iniariant divisors are different from the ones of $P_{VI}$-case.
\begin{center}
\begin{tabular}{|c||c|c|c|c|c|} \hline 
Invariant divisors &  $f_0$ & $f_1$ & $f_2$ & $f_3$ & $f_4$ \\ \hline
$P_{VI}$ &  $q-t$ & $q-\infty$ & $p$ & $q-1$ & $q$ \\ \hline
2-CPIII &  $p_1-1$ & $p_1$ & $q_1q_2-1$ & $p_2$ & $p_2-t$ \\ \hline
\end{tabular}
\end{center}

At first, let us consider $D_6^{(1)}$ case. In this paper, we present a 6-parameter family of coupled Painlev\'e III systems with affine Weyl group symmetry of type $D_6^{(1)}$ explicitly given by
\begin{equation}\label{1}
\frac{dx}{dt}=\frac{\partial H}{\partial y}, \quad \frac{dy}{dt}=-\frac{\partial H}{\partial x}, \quad \frac{dz}{dt}=\frac{\partial H}{\partial w}, \quad \frac{dw}{dt}=-\frac{\partial H}{\partial z}, \quad \frac{dq}{dt}=\frac{\partial H}{\partial p}, \quad \frac{dp}{dt}=-\frac{\partial H}{\partial q}
\end{equation}
with the polynomial Hamiltonian
\begin{align}\label{2}
\begin{split}
H &=\frac{x^2(y-1)y+x\{(\alpha_0+\alpha_1)y-\alpha_1 \}+ty}{t}\\
&+\frac{z^2(w-1)w+z\{(\alpha_0+\alpha_1+2\alpha_2+2\alpha_3)w-\alpha_3 \}+tw}{t}\\
&+\frac{q^2(p-t)p+q\{(\alpha_5+\alpha_6-1)p-t \alpha_5 \}+p}{t}+\frac{2yz(zw+\alpha_3)}{t}-\frac{2(y+w)p}{t}\\
&=H_{III}(x,y,t;\alpha_1,\alpha_0)+H_{III}(z,w,t;\alpha_3,\alpha_0+\alpha_1+2\alpha_2+\alpha_3)\\
&+{\tilde H}_{III}(q,p,t;\alpha_5,1-\alpha_6)+\frac{2yz(zw+\alpha_3)}{t}-\frac{2(y+w)p}{t}.
\end{split}
\end{align}
Here $x,y,z,w,q$ and $p$ denote unknown complex variables, and $\alpha_0,\alpha_1, \dots ,\alpha_6$ are complex parameters satisfying the relation:
\begin{equation}\label{3}
\alpha_0+\alpha_1+2\alpha_2+2\alpha_3+2\alpha_4+\alpha_5+\alpha_6=1.
\end{equation}
The symbols $H_{III},{\tilde H}_{III}$ are given by
\begin{align}\label{4,5}
&H_{III}(u,v,t;\gamma_0,\gamma_1,\gamma_2)=\frac{u^2v(v-1)+u\{(\gamma_0+\gamma_2)v-\gamma_0\}+tv}{t} \quad (\gamma_0+2\gamma_1+\gamma_2=1),\\
&{\tilde H}_{III}(U,V,t;\gamma_0,\gamma_1,\gamma_2)=\frac{U^2V(V-t)-U\{(-\gamma_0+\gamma_2)V+\gamma_0t\}+V}{t}.
\end{align}
The relation between $(u,v)$ and $(U,V)$ is given by
\begin{equation}\label{7}
(U,V)=(1/u,-u(vu+\gamma_0)).
\end{equation}
We remark that for this system we tried to seek its first integrals of polynomial type with respect to $x,y,z,w,q,p$. However, we can not find. Of course, the Hamiltonian $H$ is not the first integral.

This is the second example which gave higher order Painlev\'e type systems of type $D_{6}^{(1)}$.

We also find and study its degenerate systems with affine Weyl group symmetry of types $B_5^{(1)}$ and $D_5^{(2)}$.  In $D_5^{(2)}$-case, each differential system with respect to all principal parts has its first integral.  Nevertheless, the polynomial Hamiltonian itself is not its first integral (see Section 4).

We give an explicit confluence process from the $D_{6}^{(1)}$ system, respectively. We will show that the system of type $B_5^{(1)}$ is equivalent to the $D_5^{(1)}$ system (see \cite{Sasa6}) by explicit birational and symplectic transformations with some parameter's changes (see Section 3).

\section{The system of type $D_6^{(1)}$}
In this section, we present a 6-parameter family of coupled Painlev\'e III systems with affine Weyl group symmetry of type $D_6^{(1)}$ explicitly given by
\begin{equation}\label{4}
  \left\{
  \begin{aligned}
   \frac{dx}{dt} &=\frac{\partial H}{\partial y}=\frac{2x^2y+2z^2 w-x^2-(\alpha_0+\alpha_1)x+2\alpha_3 z-2p+t}{t},\\
   \frac{dy}{dt} &=-\frac{\partial H}{\partial x}=\frac{-2xy^2+2xy-(\alpha_0+\alpha_1)y+\alpha_1}{t},\\
   \frac{dz}{dt} &=\frac{\partial H}{\partial w}=\frac{2z^2w+2yz^2-z^2-(2\alpha_4-1+\alpha_5+\alpha_6)z-2p+t}{t},\\
   \frac{dw}{dt} &=-\frac{\partial H}{\partial z}=\frac{-2zw^2-4yzw+2zw-2\alpha_3 y+(2\alpha_4-1+\alpha_5+\alpha_6)w+\alpha_3}{t},\\
   \frac{dq}{dt} &=\frac{\partial H}{\partial p}=\frac{2q^2p-tq^2-2y-2w+(\alpha_5+\alpha_6-1)q+1}{t},\\
   \frac{dp}{dt} &=-\frac{\partial H}{\partial q}=\frac{-2qp^2+2tqp-(\alpha_5+\alpha_6-1)p+t \alpha_5}{t}
   \end{aligned}
  \right. 
\end{equation}
with the Hamiltonian \eqref{2}.

\begin{theorem}\label{2.1}
The system \eqref{4} admits extended affine Weyl group symmetry of type $D_6^{(1)}$ as the group of its B{\"a}cklund transformations whose generators are explicitly given as follows{\rm : \rm}with the notation $(*):=(x,y,z,w,q,p,t;\alpha_0,\alpha_1, \dots ,\alpha_6),$

\begin{align*}
        s_{0}: (*) &\rightarrow \left(x+\frac{\alpha_0}{y-1},y,z,w,q,p,t;-\alpha_0,\alpha_1,\alpha_2+\alpha_0,\alpha_3,\alpha_4,\alpha_5,\alpha_6 \right), \\
        s_{1}: (*) &\rightarrow \left(x+\frac{\alpha_1}{y},y,z,w,q,p,t;\alpha_0,-\alpha_1,\alpha_2+\alpha_1,\alpha_3,\alpha_4,\alpha_5,\alpha_6 \right), \\
        s_{2}: (*) &\rightarrow  \left(x,y-\frac{\alpha_2}{x-z},z,w+\frac{\alpha_2}{x-z},q,p,t;\alpha_0+\alpha_2,\alpha_1+\alpha_2,-\alpha_2,\alpha_3+\alpha_2,\alpha_4,\alpha_5,\alpha_6 \right), \\
        s_{3}: (*) &\rightarrow \left(x,y,z+\frac{\alpha_3}{w},w,q,p,t;\alpha_0,\alpha_1,\alpha_2+\alpha_3,-\alpha_3,\alpha_4+\alpha_3,\alpha_5,\alpha_6 \right),\\
        s_{4}: (*) &\rightarrow \left(x,y,z,w-\frac{\alpha_4q}{zq-1},q,p-\frac{\alpha_4z}{zq-1},t;\alpha_0,\alpha_1,\alpha_2,\alpha_3+\alpha_4,-\alpha_4,\alpha_5+\alpha_4,\alpha_6+\alpha_4 \right), \\
        s_{5}: (*) &\rightarrow \left(x,y,z,w,q+\frac{\alpha_5}{p},p,t;\alpha_0,\alpha_1,\alpha_2,\alpha_3,\alpha_4+\alpha_5,-\alpha_5,\alpha_6 \right), \\
        s_{6}: (*) &\rightarrow \left(x,y,z,w,q+\frac{\alpha_6}{p-t},p,t;\alpha_0,\alpha_1,\alpha_2,\alpha_3,\alpha_4+\alpha_6,\alpha_5,-\alpha_6 \right), \\
        \pi_1: (*) &\rightarrow (-x,1-y,-z,-w,-q,t-p,t;\alpha_1,\alpha_0,\alpha_2,\alpha_3,\alpha_4,\alpha_6,\alpha_5), \\
        \pi_2: (*) &\rightarrow \left(tq,\frac{p}{t},\frac{t}{z},-\frac{(zw+\alpha_3)z}{t},\frac{x}{t},ty,t;\alpha_6,\alpha_5,\alpha_4,\alpha_3,\alpha_2,\alpha_1,\alpha_0 \right),
        \end{align*}
        \begin{align*}
        \pi_3: (*) &\rightarrow (x,y,z,w,q,p-t,-t;\alpha_0,\alpha_1,\alpha_2,\alpha_3,\alpha_4,\alpha_6,\alpha_5), \\
        \pi_4: (*) &\rightarrow (-x,1-y,-z,-w,-q,-p,-t;\alpha_1,\alpha_0,\alpha_2,\alpha_3,\alpha_4,\alpha_5,\alpha_6).
\end{align*}
\end{theorem}
The B{\"a}cklund transformations of this system satisfy the universal description for $D_6^{(1)}$ root system. Since these universal B{\"a}cklund transformations have Lie theoretic origin, similarity reduction of a Drinfeld-Sokolov hierarchy admits such a B{\"a}cklund symmetry.

\begin{figure}
\unitlength 0.1in
\begin{picture}( 60.2400, 32.7600)( 12.3800,-39.9800)
%
\special{pn 8}%
\special{ar 1560 940 322 192  0.0000000 6.2831853}%
%
\special{pn 8}%
\special{ar 2250 1438 322 190  0.0000000 6.2831853}%
%
\special{pn 8}%
\special{ar 1630 1982 322 190  0.0000000 6.2831853}%
%
\special{pn 8}%
\special{ar 3270 1438 322 190  0.0000000 6.2831853}%
%
\special{pn 8}%
\special{ar 4270 1438 322 190  0.0000000 0.0469667}%
\special{ar 4270 1438 322 190  0.1878669 0.2348337}%
\special{ar 4270 1438 322 190  0.3757339 0.4227006}%
\special{ar 4270 1438 322 190  0.5636008 0.6105675}%
\special{ar 4270 1438 322 190  0.7514677 0.7984344}%
\special{ar 4270 1438 322 190  0.9393346 0.9863014}%
\special{ar 4270 1438 322 190  1.1272016 1.1741683}%
\special{ar 4270 1438 322 190  1.3150685 1.3620352}%
\special{ar 4270 1438 322 190  1.5029354 1.5499022}%
\special{ar 4270 1438 322 190  1.6908023 1.7377691}%
\special{ar 4270 1438 322 190  1.8786693 1.9256360}%
\special{ar 4270 1438 322 190  2.0665362 2.1135029}%
\special{ar 4270 1438 322 190  2.2544031 2.3013699}%
\special{ar 4270 1438 322 190  2.4422701 2.4892368}%
\special{ar 4270 1438 322 190  2.6301370 2.6771037}%
\special{ar 4270 1438 322 190  2.8180039 2.8649706}%
\special{ar 4270 1438 322 190  3.0058708 3.0528376}%
\special{ar 4270 1438 322 190  3.1937378 3.2407045}%
\special{ar 4270 1438 322 190  3.3816047 3.4285714}%
\special{ar 4270 1438 322 190  3.5694716 3.6164384}%
\special{ar 4270 1438 322 190  3.7573386 3.8043053}%
\special{ar 4270 1438 322 190  3.9452055 3.9921722}%
\special{ar 4270 1438 322 190  4.1330724 4.1800391}%
\special{ar 4270 1438 322 190  4.3209393 4.3679061}%
\special{ar 4270 1438 322 190  4.5088063 4.5557730}%
\special{ar 4270 1438 322 190  4.6966732 4.7436399}%
\special{ar 4270 1438 322 190  4.8845401 4.9315068}%
\special{ar 4270 1438 322 190  5.0724070 5.1193738}%
\special{ar 4270 1438 322 190  5.2602740 5.3072407}%
\special{ar 4270 1438 322 190  5.4481409 5.4951076}%
\special{ar 4270 1438 322 190  5.6360078 5.6829746}%
\special{ar 4270 1438 322 190  5.8238748 5.8708415}%
\special{ar 4270 1438 322 190  6.0117417 6.0587084}%
\special{ar 4270 1438 322 190  6.1996086 6.2465753}%
%
\special{pn 8}%
\special{pa 3590 1392}%
\special{pa 3950 1392}%
\special{dt 0.045}%
\special{sh 1}%
\special{pa 3950 1392}%
\special{pa 3884 1372}%
\special{pa 3898 1392}%
\special{pa 3884 1412}%
\special{pa 3950 1392}%
\special{fp}%
%
\special{pn 8}%
\special{pa 3580 1514}%
\special{pa 3960 1514}%
\special{dt 0.045}%
\special{sh 1}%
\special{pa 3960 1514}%
\special{pa 3894 1494}%
\special{pa 3908 1514}%
\special{pa 3894 1534}%
\special{pa 3960 1514}%
\special{fp}%
%
\special{pn 8}%
\special{pa 1790 1072}%
\special{pa 2040 1284}%
\special{fp}%
%
\special{pn 8}%
\special{pa 1800 1828}%
\special{pa 2040 1592}%
\special{fp}%
%
\special{pn 8}%
\special{pa 2570 1444}%
\special{pa 2930 1444}%
\special{fp}%
\put(14.3000,-20.2600){\makebox(0,0)[lb]{$y-1$}}%
\put(14.1000,-9.9500){\makebox(0,0)[lb]{$y$}}%
\put(20.5000,-15.0700){\makebox(0,0)[lb]{$x-z$}}%
\put(31.2000,-15.0100){\makebox(0,0)[lb]{$w$}}%
%
\special{pn 8}%
\special{ar 5952 912 322 190  0.0000000 0.0470588}%
\special{ar 5952 912 322 190  0.1882353 0.2352941}%
\special{ar 5952 912 322 190  0.3764706 0.4235294}%
\special{ar 5952 912 322 190  0.5647059 0.6117647}%
\special{ar 5952 912 322 190  0.7529412 0.8000000}%
\special{ar 5952 912 322 190  0.9411765 0.9882353}%
\special{ar 5952 912 322 190  1.1294118 1.1764706}%
\special{ar 5952 912 322 190  1.3176471 1.3647059}%
\special{ar 5952 912 322 190  1.5058824 1.5529412}%
\special{ar 5952 912 322 190  1.6941176 1.7411765}%
\special{ar 5952 912 322 190  1.8823529 1.9294118}%
\special{ar 5952 912 322 190  2.0705882 2.1176471}%
\special{ar 5952 912 322 190  2.2588235 2.3058824}%
\special{ar 5952 912 322 190  2.4470588 2.4941176}%
\special{ar 5952 912 322 190  2.6352941 2.6823529}%
\special{ar 5952 912 322 190  2.8235294 2.8705882}%
\special{ar 5952 912 322 190  3.0117647 3.0588235}%
\special{ar 5952 912 322 190  3.2000000 3.2470588}%
\special{ar 5952 912 322 190  3.3882353 3.4352941}%
\special{ar 5952 912 322 190  3.5764706 3.6235294}%
\special{ar 5952 912 322 190  3.7647059 3.8117647}%
\special{ar 5952 912 322 190  3.9529412 4.0000000}%
\special{ar 5952 912 322 190  4.1411765 4.1882353}%
\special{ar 5952 912 322 190  4.3294118 4.3764706}%
\special{ar 5952 912 322 190  4.5176471 4.5647059}%
\special{ar 5952 912 322 190  4.7058824 4.7529412}%
\special{ar 5952 912 322 190  4.8941176 4.9411765}%
\special{ar 5952 912 322 190  5.0823529 5.1294118}%
\special{ar 5952 912 322 190  5.2705882 5.3176471}%
\special{ar 5952 912 322 190  5.4588235 5.5058824}%
\special{ar 5952 912 322 190  5.6470588 5.6941176}%
\special{ar 5952 912 322 190  5.8352941 5.8823529}%
\special{ar 5952 912 322 190  6.0235294 6.0705882}%
\special{ar 5952 912 322 190  6.2117647 6.2588235}%
%
\special{pn 8}%
\special{ar 5972 1436 322 190  0.0000000 6.2831853}%
%
\special{pn 8}%
\special{ar 5972 1986 322 190  0.0000000 6.2831853}%
%
\special{pn 8}%
\special{ar 5002 1424 322 190  0.0000000 6.2831853}%
%
\special{pn 8}%
\special{ar 6942 1436 322 190  0.0000000 6.2831853}%
%
\special{pn 8}%
\special{pa 5942 1110}%
\special{pa 5942 1238}%
\special{fp}%
%
\special{pn 8}%
\special{pa 5942 1634}%
\special{pa 5942 1788}%
\special{fp}%
%
\special{pn 8}%
\special{pa 5332 1416}%
\special{pa 5642 1416}%
\special{fp}%
%
\special{pn 8}%
\special{pa 6292 1424}%
\special{pa 6602 1424}%
\special{fp}%
\put(48.5100,-14.8000){\makebox(0,0)[lb]{$w$}}%
\put(57.4100,-9.6200){\makebox(0,0)[lb]{$w-1$}}%
\put(57.0100,-14.7400){\makebox(0,0)[lb]{$zq-1$}}%
\put(57.3100,-20.4300){\makebox(0,0)[lb]{$p-t$}}%
\put(67.8100,-14.8000){\makebox(0,0)[lb]{$p$}}%
%
\special{pn 20}%
\special{pa 4330 2238}%
\special{pa 4330 2558}%
\special{fp}%
\special{sh 1}%
\special{pa 4330 2558}%
\special{pa 4350 2490}%
\special{pa 4330 2504}%
\special{pa 4310 2490}%
\special{pa 4330 2558}%
\special{fp}%
%
\special{pn 20}%
\special{pa 5210 2250}%
\special{pa 5210 2570}%
\special{fp}%
\special{sh 1}%
\special{pa 5210 2570}%
\special{pa 5230 2504}%
\special{pa 5210 2518}%
\special{pa 5190 2504}%
\special{pa 5210 2570}%
\special{fp}%
%
\special{pn 8}%
\special{ar 1582 2768 322 190  0.0000000 6.2831853}%
%
\special{pn 8}%
\special{ar 2272 3264 322 192  0.0000000 6.2831853}%
%
\special{pn 8}%
\special{ar 1652 3808 322 192  0.0000000 6.2831853}%
%
\special{pn 8}%
\special{ar 3292 3264 322 192  0.0000000 6.2831853}%
%
\special{pn 8}%
\special{pa 1812 2900}%
\special{pa 2062 3110}%
\special{fp}%
%
\special{pn 8}%
\special{pa 1822 3654}%
\special{pa 2062 3418}%
\special{fp}%
%
\special{pn 8}%
\special{pa 2592 3270}%
\special{pa 2952 3270}%
\special{fp}%
\put(14.5200,-38.5200){\makebox(0,0)[lb]{$y-1$}}%
\put(14.3200,-28.2200){\makebox(0,0)[lb]{$y$}}%
\put(20.7200,-33.3400){\makebox(0,0)[lb]{$x-z$}}%
\put(31.4200,-33.2700){\makebox(0,0)[lb]{$w$}}%
%
\special{pn 8}%
\special{ar 4270 3280 322 190  0.0000000 6.2831853}%
%
\special{pn 8}%
\special{pa 3630 3260}%
\special{pa 3940 3260}%
\special{fp}%
\put(40.0000,-33.1700){\makebox(0,0)[lb]{$zq-1$}}%
%
\special{pn 8}%
\special{ar 5180 2814 322 190  0.0000000 6.2831853}%
%
\special{pn 8}%
\special{ar 5140 3768 322 190  0.0000000 6.2831853}%
%
\special{pn 8}%
\special{pa 4480 3128}%
\special{pa 4880 2910}%
\special{fp}%
%
\special{pn 8}%
\special{pa 4520 3416}%
\special{pa 4870 3664}%
\special{fp}%
\put(49.2000,-38.2400){\makebox(0,0)[lb]{$p-t$}}%
\put(50.1000,-28.6400){\makebox(0,0)[lb]{$p$}}%
\put(17.1000,-23.5000){\makebox(0,0)[lb]{Dynkin diagram of type $B_4^{(1)}$}}%
\put(53.3000,-23.5000){\makebox(0,0)[lb]{Dynkin diagram of type $D_4^{(1)}$}}%
\put(19.6000,-41.1200){\makebox(0,0)[lb]{Dynkin diagram of type $D_6^{(1)}$}}%
\end{picture}%
\label{fig:CPIIID6fig1}
\caption{The figure denotes the Dynkin diagram of types $B_4^{(1)},D_4^{(1)}$ and $D_6^{(1)}$. The symbol in each circle denotes the invariant divisors of the systems of types $B_4^{(1)},D_4^{(1)}$ and $D_6^{(1)}$.}
\end{figure}
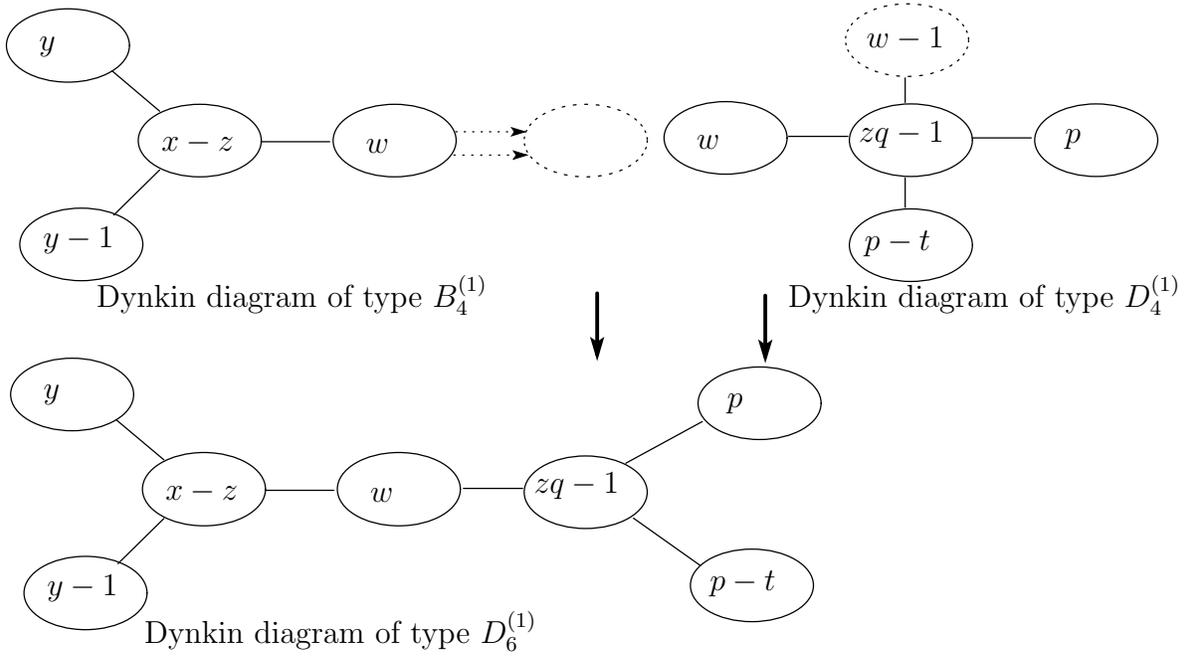

\begin{proposition}
Let us define the following translation operators{\rm : \rm}
\begin{align}
\begin{split}
&T_1:=\pi_1s_5s_4s_3s_2s_1s_0s_1s_2s_3s_4s_5, \quad T_2:=s_4s_6T_1s_6s_4,\\
&T_3:=s_6T_1s_6, \quad T_4:=\pi_2T_1\pi_2, \quad T_5:=\pi_2T_2\pi_2, \quad T_6:=\pi_2T_3\pi_2.
\end{split}
\end{align}
These translation operators act on parameters $\alpha_i$ as follows$:$
\begin{align}
\begin{split}
T_1(\alpha_0,\alpha_1,\ldots,\alpha_6)=&(\alpha_0,\alpha_1,\ldots,\alpha_6)+(0,0,0,0,0,-1,1),\\
T_2(\alpha_0,\alpha_1,\ldots,\alpha_6)=&(\alpha_0,\alpha_1,\ldots,\alpha_6)+(0,0,0,1,-1,0,0),\\
T_3(\alpha_0,\alpha_1,\ldots,\alpha_6)=&(\alpha_0,\alpha_1,\ldots,\alpha_6)+(0,0,0,0,1,-1,-1),\\
T_4(\alpha_0,\alpha_1,\ldots,\alpha_6)=&(\alpha_0,\alpha_1,\ldots,\alpha_6)+(1,-1,0,0,0,0,0),\\
T_5(\alpha_0,\alpha_1,\ldots,\alpha_6)=&(\alpha_0,\alpha_1,\ldots,\alpha_6)+(0,0,-1,1,0,0,0),\\
T_6(\alpha_0,\alpha_1,\ldots,\alpha_6)=&(\alpha_0,\alpha_1,\ldots,\alpha_6)+(-1,-1,1,0,0,0,0).
\end{split}
\end{align}
\end{proposition}

\begin{theorem}\label{2.2}
Let us consider a polynomial Hamiltonian system with Hamiltonian $H \in {\Bbb C}(t)[x,y,z,w,q,p]$. We assume that

$(A1) \ deg(H)=4$ with respect to $x,y,z,w,q,p$.

$(A2)$ This system becomes again a polynomial Hamiltonian system in each coordinate system $r_i \ (i=0,1,2,3,5,6)${\rm : \rm}
\begin{align*}
r_0&:x_0=1/x, \ y_0=-((y-1)x+\alpha_0)x, \ z_0=z, \ w_0=w, \ q_0=q, \ p_0=p, \\
r_1&:x_1=1/x, \ y_1=-(yx+\alpha_1)x, \ z_1=z, \ w_1=w, \ q_1=q, \ p_1=p, \\
r_2&:x_2=-((x-z)y-\alpha_2)y, \ y_2=1/y, \ z_2=z, \ w_2=w+y, \ q_2=q, \ p_2=p, \\
r_3&:x_3=x, \ y_3=y, \ z_3=1/z, \ w_3=-(wz+\alpha_3)z, \ q_3=q, \ p_3=p, \\
r_5&:x_5=x, \ y_5=y, \ z_5=z, \ w_5=w, \ q_5=1/q, \ p_5=-(pq+\alpha_5)q, \\
r_6&:x_6=x, \ y_6=y, \ z_6=z, \ w_6=w. \ q_6=1/q, \ p_6=-((p-t)q+\alpha_6)q.
\end{align*}

$(A3)$ In addition to the assumption $(A2)$, the Hamiltonian system in the coordinate system $r_3$ becomes again a polynomial Hamiltonian system in the coordinate system $r_4${\rm : \rm}
\begin{align*}
r_4&:x_4=x_3, \ y_4=y_3, \ z_4=-((z_3-q_3)w_3-\alpha_4)w_3, \ w_4=1/w_3, \ q_4=q_3, \ p_4=p_3+w_3.
\end{align*}
Then such a system coincides with the system \eqref{4} with the polynomial Hamiltonian \eqref{2}.
\end{theorem}
We note that the conditions $(A2)$ and $(A3)$ should be read that
\begin{align*}
&r_j(H) \quad (j=0,1,2,3,5), \quad r_6(H+q), \quad r_4(r_3(H))
\end{align*}
are polynomials with respect to $x,y,z,w,q,p$ or $x_3,y_3,z_3,w_3,q_3,p_3$.

Finally, we consider the rational and algebraic solutions of the system \eqref{4}.

At first, we consider the Dynkin diagram automorphism $\pi_1$. By this transformation, the fixed solution is derived from
\begin{align}
\begin{split}
&\alpha_0=\alpha_1, \quad \alpha_1=\alpha_0, \quad \alpha_5=\alpha_6, \quad \alpha_6=\alpha_5,\\
&x=-x, \quad y=1-y, \quad z=-z, \quad w=-w, \quad q=-q, \quad p=t-p.
\end{split}
\end{align}
Then we obtain
\begin{align}
\begin{split}
&(\alpha_0,\alpha_1,\dots,\alpha_6)=\left(\frac{1}{2}-\alpha_2-\alpha_3-\alpha_4-\alpha_6,\frac{1}{2}-\alpha_2-\alpha_3-\alpha_4-\alpha_6,\alpha_2,\alpha_3,\alpha_4,\alpha_6,\alpha_6 \right),\\
&(x,y,z,w,q,p)=\left(0,\frac{1}{2},0,0,0,\frac{t}{2} \right).
\end{split}
\end{align}

Next, we find two algebraic solutions:
\begin{align}
\begin{split}
&(\alpha_0,\alpha_1,\dots,\alpha_6)=\left(\frac{1-2\alpha_3}{2},0,0,\alpha_3,0,0,\frac{1-2\alpha_3}{2} \right),\\
&(x,y,z,w,q,p)=\left(\sqrt{t},0,\sqrt{t},-\frac{\alpha_3}{2\sqrt{t}},\frac{1}{\sqrt{t}},0 \right)
\end{split}
\end{align}
and
\begin{align}
\begin{split}
&(\alpha_0,\alpha_1,\dots,\alpha_6)=\left(0,\frac{1+2\alpha_3}{2},-\alpha_3,\alpha_3,-\alpha_3,\frac{1+2\alpha_3}{2},0 \right),\\
&(x,y,z,w,q,p)=\left(-\sqrt{t},1,\sqrt{t},-\frac{\alpha_3}{2\sqrt{t}},-\frac{1}{\sqrt{t}},t \right).
\end{split}
\end{align}

\section{The system of type $B_5^{(1)}$}

In this section, we present a 5-parameter family of coupled Painlev\'e systems with affine Weyl group symmetry of type $B_5^{(1)}$ explicitly given by
\begin{equation}\label{21}
  \left\{
  \begin{aligned}
   \frac{dx}{dt} &=\frac{\partial H}{\partial y}=\frac{2x^2y+2z^2 w+2\alpha_0 x+2\alpha_2 z-2p+t}{t},\\
   \frac{dy}{dt} &=-\frac{\partial H}{\partial x}=\frac{-2xy^2-2\alpha_0 y-1}{t},\\
   \frac{dz}{dt} &=\frac{\partial H}{\partial w}=\frac{2z^2w+2yz^2+2(\alpha_0+\alpha_1+\alpha_2)z-2p+t}{t},\\
   \frac{dw}{dt} &=-\frac{\partial H}{\partial z}=\frac{-2zw^2-4yzw-2\alpha_2 y-2(\alpha_0+\alpha_1+\alpha_2)w}{t},\\
   \frac{dq}{dt} &=\frac{\partial H}{\partial p}=\frac{2q^2p-tq^2-2y-2w-2(\alpha_0+\alpha_1+\alpha_2+\alpha_3)q}{t},\\
   \frac{dp}{dt} &=-\frac{\partial H}{\partial q}=\frac{-2qp^2+2tqp+2(\alpha_0+\alpha_1+\alpha_2+\alpha_3)p+t \alpha_4}{t}
   \end{aligned}
  \right. 
\end{equation}
with the polynomial Hamiltonian
\begin{align}\label{22}
\begin{split}
H &=\frac{x^2y^2+2\alpha_0 xy+x+ty}{t}+\frac{z^2w^2+2(\alpha_0+\alpha_1+\alpha_2)zw+tw}{t}\\
&+\frac{q^2p^2-tq^2p+(\alpha_4+\alpha_5-1)qp-\alpha_4 tq}{t}+\frac{2yz(zw+\alpha_2)}{t}-\frac{2(y+w)p}{t}\\
&=H_{III}^{D_7^{(1)}}(x,y,t;2\alpha_0)+H_1(z,w,t;2(\alpha_0+\alpha_1+\alpha_2))\\
&+H_2(q,p,t;\alpha_4+\alpha_5-1,-\alpha_4)+\frac{2yz(zw+\alpha_2)}{t}-\frac{2(y+w)p}{t}.
\end{split}
\end{align}
Here $x,y,z,w,q$ and $p$ denote unknown complex variables, and $\alpha_0,\alpha_1, \dots ,\alpha_5$ are complex parameters satisfying the relation:
\begin{equation}\label{23}
2\alpha_0+2\alpha_1+2\alpha_2+2\alpha_3+\alpha_4+\alpha_5=1.
\end{equation}
The symbols $H_{III}^{D_7^{(1)}},H_1$ and $H_2$ are given by
\begin{align}\label{24}
&H_{III}^{D_7^{(1)}}(q,p,t;\beta_1)=\frac{q^2p^2+\beta_1 qp+q+tp}{t} \quad (\beta_0+\beta_1=1),\\
&H_1(q,p,t;\alpha)=\frac{q^2p^2+\alpha qp+tp}{t},\\
&H_2(q,p,t;\alpha,\beta)=\frac{q^2p^2-tq^2p+\alpha qp+\beta tq}{t}.
\end{align}
We remark that for $y=q/{\tau}, \ t={\tau}^2$  the Hamiltonian system with $H_{III}^{D_7^{(1)}}$ is the special case of the third Painlev\'e system (see \cite{T}):
\begin{equation}
\frac{d^2y}{d{\tau}^2}=\frac{1}{y}\left(\frac{dy}{d{\tau}}\right)^2-\frac{1}{\tau}\frac{dy}{d{\tau}}+\frac{1}{\tau}(ay^2+b)+cy^3+\frac{d}{y}
\end{equation}
with
\begin{equation}
a=-8, \quad b=4(1-\beta_1), \quad c=0, \quad d=-4.
\end{equation}

From the viewpoint of symmetry, the Hamiltonian system
\begin{equation}
\frac{dq}{dt}=\frac{\partial H_{III}^{D_7^{(1)}}}{\partial p}, \quad \frac{dp}{dt}=-\frac{\partial H_{III}^{D_7^{(1)}}}{\partial q}
\end{equation}
 has extended affine Weyl group symmetry of type $A_1^{(1)}$, whose generators $<s_0,s_1,\pi=\sigma \circ s_1>$ are explicitly given as follows (see \cite{T}):
\begin{equation}
  \left\{
  \begin{aligned}
   s_0(q,p,t;\beta_0,\beta_1) &=\left(q,p+\frac{\beta_0}{q}-\frac{t}{q^2},-t;-\beta_0,\beta_1+2\beta_0 \right),\\
   s_1(q,p,t;\beta_0,\beta_1) &=\left(-q+\frac{\beta_1}{p}+\frac{1}{p^2},-p,-t;\beta_0+2\beta_1,-\beta_1 \right),\\
   \sigma(q,p,t;\beta_0,\beta_1) &=\left(tp,-\frac{q}{t},-t;\beta_1,\beta_0 \right).
   \end{aligned}
  \right. 
\end{equation}

\begin{proposition}\label{pro:1}
By the following birational and symplectic transformations $tr_i \ (i=1,2)${\rm : \rm}
\begin{equation}\label{0.1}
  \left\{
  \begin{aligned}
   tr_1(q,p) &=(q/t,tp),\\
   tr_2(q,p) &=(-p/t,tq),
   \end{aligned}
  \right. 
\end{equation}
the Hamiltonians $H_1$ and $H_2$ satisfy the following relations{\rm : \rm}
\begin{equation}
tr_1(H_1)=\frac{q^2p^2+(\alpha-1)qp+p}{t}, \quad tr_2(H_2)=\frac{q^2p^2+qp^2-(\alpha+1)qp+\beta p}{t}.
\end{equation}
\end{proposition}
Here, for notational convenience, we use the same symbol $q,p,\alpha$.

By Proposition \ref{pro:1}, we see that the Hamiltonian system with
$$
K_1:=\frac{q^2p^2+(\alpha-1)qp+p}{t}
$$
has the first integral.
\begin{proposition}
The system with the Hamiltonian $K_1$
\begin{equation}
\frac{dq}{dt}=\frac{\partial K_1}{\partial p}, \quad \frac{dp}{dt}=-\frac{\partial K_1}{\partial q}
\end{equation}
has the first integral $I_1${\rm : \rm}
\begin{equation}
I_1=q^2p^2+(\alpha-1)qp+p.
\end{equation}
\end{proposition}
We see that the relation between the Hamiltonian $K_1$ and the first integral $I_1$ is explicitly given by
\begin{equation}
tK_1=I_1.
\end{equation}
We also show that the Hamiltonian system with
$$
K_2:=\frac{q^2p^2+qp^2-(\alpha+1)qp+\beta p}{t}
$$
has the first integral.
\begin{proposition}
The system with the Hamiltonian $K_2$
\begin{equation}
\frac{dq}{dt}=\frac{\partial K_2}{\partial p}, \quad \frac{dp}{dt}=-\frac{\partial K_2}{\partial q}
\end{equation}
has the first integral $I_2${\rm : \rm}
\begin{equation}
I_2=q^2p^2+qp^2-(\alpha+1)qp+\beta p.
\end{equation}
\end{proposition}
We see that the relation between the Hamiltonian $K_2$ and the first integral $I_2$ is explicitly given by
\begin{equation}
tK_2=I_2.
\end{equation}

We remark that for this system we tried to seek its first integrals of polynomial type with respect to $x,y,z,w,q,p$. However, we can not find. Of course, the Hamiltonian $H$ is not the first integral.

\begin{theorem}\label{3.1}
The system \eqref{21} admits extended affine Weyl group symmetry of type $B_5^{(1)}$ as the group of its B{\"a}cklund transformations, whose generators are explicitly given as follows{\rm : \rm}with the notation $(*):=(x,y,z,w,q,p,t;\alpha_0,\alpha_1, \dots ,\alpha_5),$
\begin{align*}
        s_{0}: (*) &\rightarrow \left(-x-\frac{2\alpha_0}{y}-\frac{1}{y^2},-y,-z,-w,-q,-p,-t;-\alpha_0,\alpha_1+2\alpha_0,\alpha_2,\alpha_3,\alpha_4,\alpha_5 \right), \\
        s_{1}: (*) &\rightarrow \left(x,y-\frac{\alpha_1}{x-z},z,w+\frac{\alpha_1}{x-z},q,p,t;\alpha_0+\alpha_1,-\alpha_1,\alpha_2+\alpha_1,\alpha_3,\alpha_4,\alpha_5 \right), \\
        s_{2}: (*) &\rightarrow  \left(x,y,z+\frac{\alpha_2}{w},w,q,p,t;\alpha_0,\alpha_1+\alpha_2,-\alpha_2,\alpha_3+\alpha_2,\alpha_4,\alpha_5 \right),\\
        s_{3}: (*) &\rightarrow \left(x,y,z,w-\frac{\alpha_3q}{zq-1},q,p-\frac{\alpha_3z}{zq-1},t;\alpha_0,\alpha_1,\alpha_2+\alpha_3,-\alpha_3,\alpha_4+\alpha_3,\alpha_5+\alpha_3 \right), \\
        s_{4}: (*) &\rightarrow \left(x,y,z,w,q+\frac{\alpha_4}{p},p,t;\alpha_0,\alpha_1,\alpha_2,\alpha_3+\alpha_4,-\alpha_4,\alpha_5 \right), \\
        s_{5}: (*) &\rightarrow \left(x,y,z,w,q+\frac{\alpha_5}{p-t},p,t;\alpha_0,\alpha_1,\alpha_2,\alpha_3+\alpha_5,\alpha_4,-\alpha_5 \right), \\
        \pi: (*) &\rightarrow (x,y,z,w,q,p-t,-t;\alpha_0,\alpha_1,\alpha_2,\alpha_3,\alpha_5,\alpha_4).
\end{align*}
\end{theorem}

\begin{proposition}
Let us define the following translation operators{\rm : \rm}
\begin{align}
\begin{split}
&T_1:=\pi s_4 s_3 s_2 s_1 s_0 s_1 s_2 s_3 s_4, \quad T_2:=\pi s_5s_4s_3s_2s_1s_0s_1s_2s_3,\\
&T_3:=s_3 s_5 T_1 s_5 s_3, \quad T_4:=s_2 T_3 s_2, \quad T_5:=s_1 T_4 s_1.
\end{split}
\end{align}
These translation operators act on parameters $\alpha_i$ as follows$:$
\begin{align}
\begin{split}
T_1(\alpha_0,\alpha_1,\ldots,\alpha_5)=&(\alpha_0,\alpha_1,\ldots,\alpha_5)+(0,0,0,0,-1,1),\\
T_2(\alpha_0,\alpha_1,\ldots,\alpha_5)=&(\alpha_0,\alpha_1,\ldots,\alpha_5)+(0,0,0,-1,1,1),\\
T_3(\alpha_0,\alpha_1,\ldots,\alpha_5)=&(\alpha_0,\alpha_1,\ldots,\alpha_5)+(0,0,1,-1,0,0),\\
T_4(\alpha_0,\alpha_1,\ldots,\alpha_5)=&(\alpha_0,\alpha_1,\ldots,\alpha_5)+(0,1,-1,0,0,0),\\
T_5(\alpha_0,\alpha_1,\ldots,\alpha_5)=&(\alpha_0,\alpha_1,\ldots,\alpha_5)+(1,-1,0,0,0,0).
\end{split}
\end{align}
\end{proposition}

\begin{theorem}\label{3.2}
Let us consider a polynomial Hamiltonian system with Hamiltonian $H \in {\Bbb C}(t)[x,y,z,w,q,p]$. We assume that

$(B1)$ $deg(H)=4$ with respect to $x,y,z,w,q,p$.

$(B2)$ This system becomes again a polynomial Hamiltonian system in each coordinate system $r_i \ (i=0,1,2,4,5)${\rm : \rm}
\begin{align*}
r_0&:x_0=x+\frac{2\alpha_0}{y}+\frac{1}{y^2}, \ y_0=y, \ z_0=z, \ w_0=w, \ q_0=q, \ p_0=p, \\
r_1&:x_1=-((x-z)y-\alpha_1)y, \ y_1=1/y, \ z_1=z, \ w_1=w+y, \ q_1=q, \ p_1=p, \\
r_2&:x_2=x, \ y_2=y, \ z_2=1/z, \ w_2=-(wz+\alpha_2)z, \ q_2=q, \ p_2=p, \\
r_4&:x_4=x, \ y_4=y, \ z_4=z, \ w_4=w, \ q_4=1/q, \ p_4=-(pq+\alpha_4)q, \\
r_5&:x_5=x, \ y_5=y, \ z_5=z, \ w_5=w. \ q_5=1/q, \ p_5=-((p-t)q+\alpha_5)q.
\end{align*}

$(B3)$ In addition to the assumption $(B2)$, the Hamiltonian system in the coordinate system $(x_4,y_4,z_4,w_4,q_4,p_4)$ becomes again a polynomial Hamiltonian system in the coordinate system $r_3$\rm{:\rm}
\begin{align*}
r_3:x_3=x_4, \ y_3=y_4, \ z_3=-((z_4-q_4)w_4-\alpha_3)w_4, \ w_3=1/w_4, \ q_3=q_4, \ p_3=p_4+w_4.
\end{align*}
Then such a system coincides with the system \eqref{21} with the polynomial Hamiltonian \eqref{22}.
\end{theorem}
We note that the conditions $(B2)$ and $(B3)$ should be read that
\begin{align*}
&r_j(H) \quad (j=0,1,2,4), \quad r_5(H+q), \quad r_3(r_4(H))
\end{align*}
are polynomials with respect to $x,y,z,w,q,p$ or $x_4,y_4,z_4,w_4,q_4,p_4$.

Theorems \ref{3.1} and \ref{3.2} can be checked by a direct calculation, respectively.

Next, we show the confluence process from the system \eqref{4} to the system \eqref{21}.
\begin{theorem}\label{3.3}
For the system \eqref{4} of type $D_6^{(1)}$, we make the change of parameters and variables
\begin{gather}
\begin{gathered}\label{28}
\alpha_0=\frac{1}{\varepsilon}+2A_0, \quad \alpha_1=-\frac{1}{\varepsilon}, \quad \alpha_2=A_1, \quad \alpha_3=A_2, \quad \alpha_4=A_3, \quad \alpha_5=A_4, \quad \alpha_6=A_5,\\
\end{gathered}\\
\begin{gathered}\label{29}
t=\varepsilon T, \quad x=\varepsilon X,\quad y=\frac{Y}{\varepsilon}, \quad z=\varepsilon Z, \quad w=\frac{W}{\varepsilon}, \quad q=\frac{Q}{\varepsilon}, \quad p=\varepsilon P
\end{gathered}
\end{gather}
from $\alpha_0,\alpha_1, \dots ,\alpha_6,x,y,z,w,q,p$ to $A_0,A_1,\dots ,A_5,X,Y,Z,W,Q,P$. Then the system \eqref{4} can also be written in the new variables $X,Y,Z,W,Q,P$ and parameters $A_0,A_1,\dots ,A_5$ as a Hamiltonian system. This new system tends to the system \eqref{21} with the Hamiltonian \eqref{22} as $\varepsilon \rightarrow 0$.
\end{theorem}

Finally, we show the relation between the system \eqref{21} and the system of type $D_5^{(1)}$ (see \cite{Sasa6}).
\begin{theorem}\label{3.4}
For the system \eqref{21} of type $B_5^{(1)}$, we make the change of parameters and variables
\begin{gather}
\begin{gathered}\label{221}
\alpha_0=\frac{A_0-A_1}{2}, \quad \alpha_1=A_1, \quad \alpha_2=A_2, \quad \alpha_3=A_3, \quad \alpha_4=A_4, \quad \alpha_5=A_5,\\
\end{gathered}\\
\begin{gathered}\label{222}
X=((x-z)y-\alpha_1)y, \quad Y=-1/y, \quad Z=z, \quad W=w+y, \quad Q=q, \quad P=p
\end{gathered}
\end{gather}
from $\alpha_0,\alpha_1, \dots ,\alpha_5,x,y,z,w,q,p$ to $A_0,A_1,\dots ,A_5,X,Y,Z,W,Q,P$. Then the system \eqref{4} can also be written in the new variables $X,Y,Z,W,Q,P$ and parameters $A_0,A_1,\dots ,A_5$ as a Hamiltonian system. This new system tends to
\begin{equation}
\frac{dx}{dt}=\frac{\partial H}{\partial y}, \quad \frac{dy}{dt}=-\frac{\partial H}{\partial x}, \quad \frac{dz}{dt}=\frac{\partial H}{\partial w}, \quad \frac{dw}{dt}=-\frac{\partial H}{\partial z}, \quad \frac{dq}{dt}=\frac{\partial H}{\partial p}, \quad \frac{dp}{dt}=-\frac{\partial H}{\partial q}
\end{equation}
with the polynomial Hamiltonian
\begin{align}
\begin{split}
H =&\frac{x^2y^2+xy^2-(\alpha_0+\alpha_1)xy-\alpha_0y}{t}+\frac{z^2w^2+(\alpha_0+\alpha_1+2\alpha_2)zw+z+tw}{t}\\
&+\frac{q^2p^2-tq^2p-(1-\alpha_4-\alpha_5)qp-\alpha_4tq}{t}+\frac{2(xz-wp)}{t}\\
&(\alpha_0+\alpha_1+2\alpha_2+2\alpha_3+\alpha_4+\alpha_5=1).
\end{split}
\end{align}
\end{theorem}
Here, for notational convenience, we have renamed $A_i,X,Y,Z,W,Q,P$ to $\alpha_i,x,y,z,w,q,p$ (which are not the same as the previous $\alpha_i,x,y,z,w,q,p$).

We recall that this system admits extended affine Weyl group symmetry of type $D_5^{(1)}$ (see \cite{Sasa6}) as the group of its B{\"a}cklund transformations whose generators are explicitly given as follows{\rm : \rm}with the notation $(*):=(x,y,z,w,q,p,t;\alpha_0,\alpha_1, \dots ,\alpha_5),$
\begin{figure}[ht]
\unitlength 0.1in
\begin{picture}(51.30,21.92)(13.00,-25.52)
%
\special{pn 20}%
\special{ar 2462 834 244 244  0.0000000 6.2831853}%
%
\special{pn 20}%
\special{ar 2473 2121 244 244  0.0000000 6.2831853}%
%
\special{pn 20}%
\special{ar 3397 1494 244 244  0.0000000 6.2831853}%
%
\special{pn 20}%
\special{pa 2660 988}%
\special{pa 3210 1307}%
\special{fp}%
%
\special{pn 20}%
\special{pa 2638 1934}%
\special{pa 3221 1681}%
\special{fp}%
\put(23.3000,-22.3000){\makebox(0,0)[lb]{$x$}}%
\put(22.6000,-9.3000){\makebox(0,0)[lb]{$x+1$}}%
\put(31.7000,-15.8000){\makebox(0,0)[lb]{$yw+1$}}%
\put(24.2000,-18.6000){\makebox(0,0)[lb]{$0$}}%
\put(24.2000,-5.3000){\makebox(0,0)[lb]{$1$}}%
\put(32.8000,-12.1000){\makebox(0,0)[lb]{$2$}}%
%
\special{pn 20}%
\special{pa 3640 1480}%
\special{pa 4300 1480}%
\special{fp}%
%
\special{pn 20}%
\special{ar 4570 1490 244 244  0.0000000 6.2831853}%
%
\special{pn 20}%
\special{pa 4750 1330}%
\special{pa 5520 980}%
\special{fp}%
%
\special{pn 20}%
\special{pa 4720 1670}%
\special{pa 5520 1950}%
\special{fp}%
%
\special{pn 20}%
\special{ar 5740 850 244 244  0.0000000 6.2831853}%
%
\special{pn 20}%
\special{ar 5710 2110 244 244  0.0000000 6.2831853}%
\put(43.3000,-15.9000){\makebox(0,0)[lb]{$zq-1$}}%
\put(55.1000,-22.0000){\makebox(0,0)[lb]{$p-t$}}%
\put(56.3000,-9.6000){\makebox(0,0)[lb]{$p$}}%
\put(44.5000,-12.2000){\makebox(0,0)[lb]{$3$}}%
\put(57.6000,-18.6000){\makebox(0,0)[lb]{$5$}}%
\put(57.6000,-5.9000){\makebox(0,0)[lb]{$4$}}%
%
\special{pn 8}%
\special{pa 1490 1510}%
\special{pa 3120 1510}%
\special{dt 0.045}%
\special{pa 3120 1510}%
\special{pa 3119 1510}%
\special{dt 0.045}%
%
\special{pn 8}%
\special{pa 4810 1480}%
\special{pa 6350 1480}%
\special{dt 0.045}%
\special{pa 6350 1480}%
\special{pa 6349 1480}%
\special{dt 0.045}%
%
\special{pn 20}%
\special{pa 1566 1408}%
\special{pa 1550 1377}%
\special{pa 1532 1349}%
\special{pa 1510 1327}%
\special{pa 1484 1313}%
\special{pa 1452 1310}%
\special{pa 1418 1315}%
\special{pa 1386 1327}%
\special{pa 1357 1344}%
\special{pa 1335 1366}%
\special{pa 1319 1392}%
\special{pa 1307 1422}%
\special{pa 1301 1454}%
\special{pa 1299 1488}%
\special{pa 1302 1523}%
\special{pa 1308 1560}%
\special{pa 1318 1596}%
\special{pa 1331 1633}%
\special{pa 1346 1668}%
\special{pa 1364 1701}%
\special{pa 1385 1728}%
\special{pa 1410 1745}%
\special{pa 1441 1747}%
\special{pa 1474 1734}%
\special{pa 1502 1711}%
\special{pa 1516 1681}%
\special{pa 1516 1668}%
\special{sp}%
%
\special{pn 20}%
\special{pa 6070 1440}%
\special{pa 6090 1406}%
\special{pa 6111 1374}%
\special{pa 6132 1344}%
\special{pa 6153 1319}%
\special{pa 6176 1299}%
\special{pa 6199 1285}%
\special{pa 6224 1280}%
\special{pa 6250 1283}%
\special{pa 6276 1295}%
\special{pa 6303 1313}%
\special{pa 6328 1336}%
\special{pa 6353 1364}%
\special{pa 6375 1395}%
\special{pa 6394 1427}%
\special{pa 6409 1461}%
\special{pa 6421 1494}%
\special{pa 6428 1527}%
\special{pa 6430 1560}%
\special{pa 6428 1592}%
\special{pa 6420 1623}%
\special{pa 6407 1654}%
\special{pa 6388 1683}%
\special{pa 6363 1710}%
\special{pa 6335 1731}%
\special{pa 6307 1740}%
\special{pa 6281 1734}%
\special{pa 6258 1714}%
\special{pa 6237 1686}%
\special{pa 6220 1660}%
\special{sp}%
%
\special{pn 20}%
\special{pa 6230 1680}%
\special{pa 6180 1600}%
\special{fp}%
\special{sh 1}%
\special{pa 6180 1600}%
\special{pa 6198 1667}%
\special{pa 6208 1645}%
\special{pa 6232 1646}%
\special{pa 6180 1600}%
\special{fp}%
\put(13.2000,-12.9000){\makebox(0,0)[lb]{$\pi_1$}}%
\put(61.5000,-12.6000){\makebox(0,0)[lb]{$\pi_2$}}%
%
\special{pn 20}%
\special{pa 1520 1680}%
\special{pa 1600 1560}%
\special{fp}%
\special{sh 1}%
\special{pa 1600 1560}%
\special{pa 1546 1604}%
\special{pa 1570 1604}%
\special{pa 1580 1627}%
\special{pa 1600 1560}%
\special{fp}%
%
\special{pn 8}%
\special{pa 3980 390}%
\special{pa 3980 2510}%
\special{dt 0.045}%
\special{pa 3980 2510}%
\special{pa 3980 2509}%
\special{dt 0.045}%
\put(43.2000,-26.8000){\makebox(0,0)[lb]{$\pi_3$}}%
%
\special{pn 20}%
\special{pa 3870 2190}%
\special{pa 3835 2200}%
\special{pa 3800 2210}%
\special{pa 3768 2223}%
\special{pa 3739 2238}%
\special{pa 3715 2256}%
\special{pa 3696 2279}%
\special{pa 3684 2306}%
\special{pa 3679 2338}%
\special{pa 3682 2372}%
\special{pa 3692 2404}%
\special{pa 3710 2433}%
\special{pa 3734 2460}%
\special{pa 3765 2483}%
\special{pa 3800 2502}%
\special{pa 3839 2519}%
\special{pa 3880 2532}%
\special{pa 3924 2542}%
\special{pa 3968 2548}%
\special{pa 4012 2551}%
\special{pa 4054 2551}%
\special{pa 4094 2548}%
\special{pa 4131 2541}%
\special{pa 4164 2531}%
\special{pa 4191 2517}%
\special{pa 4212 2500}%
\special{pa 4225 2480}%
\special{pa 4230 2456}%
\special{pa 4227 2429}%
\special{pa 4216 2400}%
\special{pa 4201 2369}%
\special{pa 4184 2337}%
\special{pa 4180 2330}%
\special{sp}%
%
\special{pn 20}%
\special{pa 4170 2320}%
\special{pa 4140 2270}%
\special{fp}%
\special{sh 1}%
\special{pa 4140 2270}%
\special{pa 4157 2337}%
\special{pa 4167 2316}%
\special{pa 4191 2317}%
\special{pa 4140 2270}%
\special{fp}%
\end{picture}%
\label{fig:D52}
\caption{Dynkin diagram of type $D_5^{(1)}$}
\end{figure}
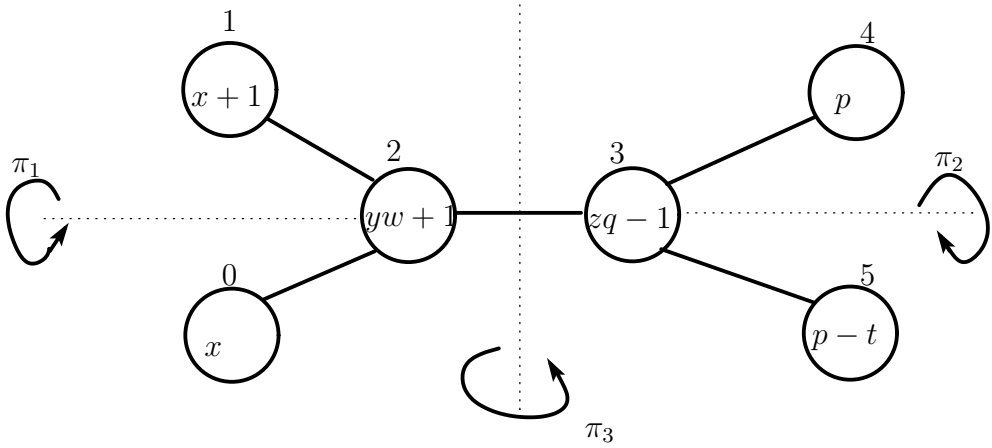
\begin{align*}
        s_{0}: (*) &\rightarrow \left(x,y-\frac{\alpha_0}{x},z,w,q,p,t;-\alpha_0,\alpha_1,\alpha_2+\alpha_0,\alpha_3,\alpha_4,\alpha_5 \right), \\
        s_{1}: (*) &\rightarrow \left(x,y-\frac{\alpha_1}{x+1},z,w,q,p,t;\alpha_0,-\alpha_1,\alpha_2+\alpha_1,\alpha_3,\alpha_4,\alpha_5 \right),\\
        s_{2}: (*) &\rightarrow  \left(x+\frac{\alpha_2w}{yw+1},y,z+\frac{\alpha_2y}{yw+1},w,q,p,t;\alpha_0+\alpha_2,\alpha_1+\alpha_2,-\alpha_2,\alpha_3+\alpha_2,\alpha_4,\alpha_5 \right),\\
        s_{3}: (*) &\rightarrow \left(x,y,z,w-\frac{\alpha_3q}{zq-1},q,p-\frac{\alpha_3z}{zq-1},t;\alpha_0,\alpha_1,\alpha_2+\alpha_3,-\alpha_3,\alpha_4+\alpha_3,\alpha_5+\alpha_3 \right), \\
        s_{4}: (*) &\rightarrow \left(x,y,z,w,q+\frac{\alpha_4}{p},p,t;\alpha_0,\alpha_1,\alpha_2,\alpha_3+\alpha_4,-\alpha_4,\alpha_5 \right), \\
        s_{5}: (*) &\rightarrow \left(x,y,z,w,q+\frac{\alpha_5}{p-t},p,t;\alpha_0,\alpha_1,\alpha_2,\alpha_3+\alpha_5,\alpha_4,-\alpha_5 \right), 
        \end{align*}
        \begin{align*}
        \pi_1: (*) &\rightarrow (-x-1,-y,-z,-w,-q,-p,-t;\alpha_1,\alpha_0,\alpha_2,\alpha_3,\alpha_4,\alpha_5), \\
        \pi_2: (*) &\rightarrow (x,y,z,w,q,p-t,-t;\alpha_0,\alpha_1,\alpha_2,\alpha_3,\alpha_5,\alpha_4), \\
        \pi_3: (*) &\rightarrow \left(\frac{(p-t)}{t},-tq,-tw,\frac{z}{t},\frac{y}{t},-t(x+1),-t;\alpha_5,\alpha_4,\alpha_3,\alpha_2,\alpha_1,\alpha_0 \right).
\end{align*}

\section{The system of type $D_5^{(2)}$ }

In this section, we present a 4-parameter family of coupled Hamiltonian systems in dimension six with extended affine Weyl group symmetry of type $D_5^{(2)}$ given by
\begin{equation}\label{41}
  \left\{
  \begin{aligned}
   \frac{dx}{dt} &=\frac{\partial H}{\partial y}=\frac{x^2y+z^2 w+\alpha_0 x+\alpha_2 z-p}{t},\\
   \frac{dy}{dt} &=-\frac{\partial H}{\partial x}=\frac{-2xy^2-2\alpha_0 y-1}{2t},\\
   \frac{dz}{dt} &=\frac{\partial H}{\partial w}=\frac{z^2w+yz^2+(\alpha_0+\alpha_1+\alpha_2)z-p}{t},\\
   \frac{dw}{dt} &=-\frac{\partial H}{\partial z}=\frac{-zw^2-2yzw-\alpha_2 y-(\alpha_0+\alpha_1+\alpha_2)w}{t},\\
   \frac{dq}{dt} &=\frac{\partial H}{\partial p}=\frac{q^2p-y-w+(\alpha_4-1)q}{t},\\
   \frac{dp}{dt} &=-\frac{\partial H}{\partial q}=\frac{-2qp^2-2(\alpha_4-1)p+t}{2t}
   \end{aligned}
  \right. 
\end{equation}
with the polynomial Hamiltonian
\begin{align}\label{42}
\begin{split}
H &=\frac{x^2y^2+2\alpha_0 xy+x}{2t}+\frac{z^2w^2+2(\alpha_0+\alpha_1+\alpha_2)zw}{2t}\\
&+\frac{q^2p^2+2(\alpha_4-1)qp-tq}{2t}+\frac{yz(zw+\alpha_2)}{t}-\frac{(y+w)p}{t}\\
&=H_3(x,y,t;2\alpha_0)+H_4(z,w,t;2(\alpha_0+\alpha_1+\alpha_2))\\
&+H_5(q,p,t;2(\alpha_4-1))+\frac{yz(zw+\alpha_2)}{t}-\frac{(y+w)p}{t}.
\end{split}
\end{align}
Here $x,y,z,w,q$ and $p$ denote unknown complex variables and $\alpha_0,\alpha_1, \dots ,\alpha_4$ are complex parameters satisfying the relation:
\begin{equation}\label{43}
\alpha_0+\alpha_1+\alpha_2+\alpha_3+\alpha_4=1.
\end{equation}
The symbols $H_3,H_4$ and $H_5$ are given by
\begin{align}\label{46}
&H_3(q,p,t;\alpha)=\frac{q^2p^2+\alpha qp+q}{2t},\\
&H_4(q,p,t;\alpha)=\frac{q^2p^2+\alpha qp}{2t},\\
&H_5(q,p,t;\alpha,\beta)=\frac{q^2p^2+\alpha qp-tq}{2t}.
\end{align}

\begin{proposition}
The system with the Hamiltonian $H_3$
\begin{equation}
\frac{dq}{dt}=\frac{\partial H_3}{\partial p}, \quad \frac{dp}{dt}=-\frac{\partial H_3}{\partial q}
\end{equation}
has the first integral $I_3${\rm : \rm}
\begin{equation}
I_3=q^2p^2+\alpha qp+q.
\end{equation}
\end{proposition}

\begin{proposition}
The system with the Hamiltonian $H_4$
\begin{equation}
\frac{dq}{dt}=\frac{\partial H_4}{\partial p}, \quad \frac{dp}{dt}=-\frac{\partial H_4}{\partial q}
\end{equation}
has the first integral $I_4${\rm : \rm}
\begin{equation}
I_4=qp.
\end{equation}
\end{proposition}

\begin{proposition}\label{pro:4.1}
By the following birational and symplectic transformation $tr_5${\rm : \rm}
\begin{equation}
   tr_5(q,p)=(tq,p/t),
\end{equation}
the Hamiltonians $H_5$ satisfy the following relation{\rm : \rm}
\begin{equation}
tr_5(H_5)=\frac{q^2p^2+\alpha qp-q}{2t}.
\end{equation}
\end{proposition}

By Proposition \ref{pro:4.1}, we see that the Hamiltonian system with
$$
K_5:=\frac{q^2p^2+\alpha qp-q}{2t}
$$
has the first integral.
\begin{proposition}
The system with the Hamiltonian $K_5$
\begin{equation}
\frac{dq}{dt}=\frac{\partial K_5}{\partial p}, \quad \frac{dp}{dt}=-\frac{\partial K_5}{\partial q}
\end{equation}
has the first integral $I_5${\rm : \rm}
\begin{equation}
I_5=q^2p^2+\alpha qp-q.
\end{equation}
\end{proposition}

In this case, each differential system with respect to all principal parts $H_3,H_4$ and $H_5$ has its first integral.  Nevertheless, the Hamiltonian $H$ is not the first integral.
For this system we tried to seek its first integrals of polynomial type with respect to $x,y,z,w,q,p$. However, we can not find.

\begin{theorem}\label{4.1}
The system \eqref{41} admits extended affine Weyl group symmetry of type $D_5^{(2)}$ as the group of its B{\"a}cklund transformations, whose generators are explicitly given as follows{\rm : \rm}with the notation $(*):=(x,y,z,w,q,p,t;\alpha_0,\alpha_1, \dots ,\alpha_4),$
\begin{align*}
        s_{0}: (*) &\rightarrow \left(-x-\frac{2\alpha_0}{y}-\frac{1}{y^2},-y,-z,-w,-q,-p,-t;-\alpha_0,\alpha_1+2\alpha_0,\alpha_2,\alpha_3,\alpha_4 \right), \\
        s_{1}: (*) &\rightarrow \left(x,y-\frac{\alpha_1}{x-z},z,w+\frac{\alpha_1}{x-z},q,p,t;\alpha_0+\alpha_1,-\alpha_1,\alpha_2+\alpha_1,\alpha_3,\alpha_4 \right), \\
        s_{2}: (*) &\rightarrow  \left(x,y,z+\frac{\alpha_2}{w},w,q,p,t;\alpha_0,\alpha_1+\alpha_2,-\alpha_2,\alpha_3+\alpha_2,\alpha_4 \right), \\
        s_{3}: (*) &\rightarrow \left(x,y,z,w-\frac{\alpha_3 q}{zq-1},q,p-\frac{\alpha_3 z}{zq-1},t;\alpha_0,\alpha_1,\alpha_2+\alpha_3,-\alpha_3,\alpha_4+\alpha_3 \right),\\
        s_{4}: (*) &\rightarrow \left(x,y,z,w,q+\frac{2\alpha_4}{p}-\frac{t}{p^2},p,-t;\alpha_0,\alpha_1,\alpha_2,\alpha_3+2\alpha_4,-\alpha_4 \right),\\
        \pi: (*) &\rightarrow \left(-tq,-\frac{p}{t},-\frac{t}{z},\frac{(zw+\alpha_2)z}{t},-\frac{x}{t},-ty,t;\alpha_4,\alpha_3,\alpha_2,\alpha_1,\alpha_0 \right).
\end{align*}
\end{theorem}

\begin{proposition}
Let us define the following translation operators{\rm : \rm}
\begin{align}
\begin{split}
&T_1:=s_4 s_3 s_2 s_1 s_0 s_1 s_2 s_3, \quad T_2:=s_3 T_1 s_3,\\
&T_3:=s_2 T_2 s_2, \quad T_4:=s_1 T_3 s_1.
\end{split}
\end{align}
These translation operators act on parameters $\alpha_i$ as follows$:$
\begin{align}
\begin{split}
T_1(\alpha_0,\alpha_1,\ldots,\alpha_4)=&(\alpha_0,\alpha_1,\ldots,\alpha_4)+(0,0,0,-2,2),\\
T_2(\alpha_0,\alpha_1,\ldots,\alpha_4)=&(\alpha_0,\alpha_1,\ldots,\alpha_4)+(0,0,-2,2,0),\\
T_3(\alpha_0,\alpha_1,\ldots,\alpha_4)=&(\alpha_0,\alpha_1,\ldots,\alpha_4)+(0,-2,2,0,0),\\
T_4(\alpha_0,\alpha_1,\ldots,\alpha_4)=&(\alpha_0,\alpha_1,\ldots,\alpha_4)+(-2,2,0,0,0).
\end{split}
\end{align}
\end{proposition}

\begin{theorem}\label{4.2}
Let us consider a polynomial Hamiltonian system with Hamiltonian $H \in {\Bbb C}(t)[x,y,z,w,q,p]$. We assume that

$(C1)$ $deg(H)=4$ with respect to $x,y,z,w,q,p$.

$(C2)$ This system becomes again a polynomial Hamiltonian system in each coordinate system $r_i \ (i=0,1,2,4)${\rm : \rm}
\begin{align*}
r_0&:x_0=x+\frac{2\alpha_0}{y}+\frac{1}{y^2}, \ y_0=y, \ z_0=z, \ w_0=w, \ q_0=q, \ p_0=p, \\
r_1&:x_1=-((x-z)y-\alpha_1)y, \ y_1=1/y, \ z_1=z, \ w_1=w+y, \ q_1=q, \ p_1=p, \\
r_2&:x_2=x, \ y_2=y, \ z_2=1/z, \ w_2=-(wz+\alpha_2)z, \ q_2=q, \ p_2=p, \\
r_4&:x_4=x, \ y_4=y, \ z_4=z, \ w_4=w. \ q_4=q+\frac{2\alpha_4}{p}-\frac{t}{p^2}, \ p_4=p.
\end{align*}

$(C3)$ In addition to the assumption $(C2)$, the Hamiltonian system in the coordinate system $(x_2,y_2,z_2,w_2,q_2,p_2)$ becomes again a polynomial Hamiltonian system in the coordinate system $r_3$\rm{:\rm}
\begin{align*}
r_3:x_3=x_2, \ y_3=y_2, \ z_3=-((z_2-q_2)w_2-\alpha_3)w_2, \ w_3=1/w_2, \ q_3=q_2, \ p_3=p_2+w_2.
\end{align*}
Then such a system coincides with the system \eqref{41} with the polynomial Hamiltonian \eqref{42}.
\end{theorem}
We note that the conditions $(C2)$ and $(C3)$ should be read that
\begin{align*}
&r_j(H) \quad (j=0,1,2), \quad r_4(H+1/p), \quad r_3(r_2(H))
\end{align*}
are polynomials with respect to $x,y,z,w,q,p$ or $x_2,y_2,z_2,w_2,q_2,p_2$.

Theorems \ref{4.1} and \ref{4.2} can be checked by a direct calculation, respectively.

Next, we show the confluence process from the system \eqref{4} to the system \eqref{41}.
\begin{theorem}\label{4.3}
For the system \eqref{4} of type $D_6^{(1)}$, we make the change of parameters and variables
\begin{gather}
\begin{gathered}\label{47}
\alpha_0=-\frac{1}{\varepsilon}+A_0, \quad \alpha_1=\frac{1}{\varepsilon}, \quad \alpha_2=\frac{A_1}{2}, \quad \alpha_3=\frac{A_2}{2}, \quad \alpha_4=\frac{A_3}{2}, \quad \alpha_5=\frac{1}{\varepsilon}, \quad \alpha_6=-\frac{1}{\varepsilon}+A_4,\\
\end{gathered}\\
\begin{gathered}\label{48}
t=\frac{\varepsilon^2}{16} T, \quad x=\frac{\varepsilon}{4} X,\quad y=\frac{2Y}{\varepsilon}, \quad z=\frac{\varepsilon}{4} Z, \quad w=\frac{2W}{\varepsilon}, \quad q=\frac{4Q}{\varepsilon}, \quad p=\frac{\varepsilon}{8} P
\end{gathered}
\end{gather}
from $\alpha_0,\alpha_1, \dots ,\alpha_6,x,y,z,w,q,p$ to $A_0,A_1,\dots ,A_4,X,Y,Z,W,Q,P$. Then the system \eqref{4} can also be written in the new variables $X,Y,Z,W,Q,P$ and parameters $A_0,A_1,\dots ,A_4$ as a Hamiltonian system. This new system tends to the system \eqref{41} with the Hamiltonian \eqref{42} as $\varepsilon \rightarrow 0$.
\end{theorem}

Finally, we find an algebraic solution of the system \eqref{41}:

\begin{align}
\begin{split}
&(\alpha_0,\alpha_1,\dots,\alpha_4)=\left(\alpha_0,\alpha_1,1-2\alpha_0-2\alpha_1,\alpha_1,\alpha_0 \right),\\
&(x,y,z,w,q,p)=(-\left(\frac{1+\sqrt{-1}}{4} \right)(t^{\frac{1}{4}}+2(1+\sqrt{-1})\sqrt{t}),-\frac{1-\sqrt{-1}}{2t^{\frac{1}{4}}},\sqrt{-t},\\
&-\frac{\sqrt{-1}(2\alpha_0+2\alpha_1-1)}{2\sqrt{t}},\frac{(1+\sqrt{-1})+4\sqrt{-1} t^{\frac{1}{4}}}{4t^{\frac{3}{4}}},\left(\frac{1}{2}-\frac{\sqrt{-1}}{2} \right)t^{\frac{3}{4}}).
\end{split}
\end{align}


\begin{thebibliography}{99}
\bibitem[1]{1} P. Painlev\'e, {\em M\'emoire sur les \'equations diff\'erentielles dont l'int\'egrale g\'en\'erale est uniforme}, Bull. Soci\'et\'e Math\'ematique de France. {\bf 28} (1900),  201--261.

\bibitem[2]{2} P. Painlev\'e, {\em Sur les \'equations diff\'erentielles du second ordre et d'ordre sup\'erieur dont l'int\'egrale est uniforme}, Acta Math. {\bf 25} (1902), 1--85. 

\bibitem[3]{3} B. Gambier, {\em Sur les \'equations diff\'erentielles du second ordre et du premier degr\'e dont l'int\'egrale g\'en\'erale est \`a points critiques fixes}, Acta Math. {\bf 33} (1910), 1--55.


\bibitem[4]{Cosgrove1} C. M. Cosgrove and G. Scoufis,
{\em Painlev\'e classification of a class of differential equations of the second order and second degree}, Studies in Applied Mathematics. {\bf 88} (1993), 25-87.

\bibitem[5]{Cosgrove2} C. M. Cosgrove,
{\em All binomial-type Painlev\'e equations of the second order and degree three or higher}, Studies in Applied Mathematics. {\bf 90} (1993), 119-187.

\bibitem[6]{6} F. Bureau, 
{\em Integration of some nonlinear systems of ordinary differential equations}, 
Annali di Matematica. {\bf 94} (1972), 345--359. 

\bibitem[7]{7} J. Chazy, 
{\em Sur les \'equations diff\'erentielles dont l'int\'egrale g\'en\'erale est uniforme et admet des singularit\'es essentielles mobiles}, 
Comptes Rendus de l'Acad\'emie des Sciences, Paris. {\bf 149} (1909), 563--565. 

\bibitem[8]{8} J. Chazy, 
{\em Sur les \'equations diff\'erentielles dont l'int\'egrale g\'en\'erale poss\'ede une coupure essentielle mobile }, 
Comptes Rendus de l'Acad\'emie des Sciences, Paris. {\bf 150} (1910), 456--458. 


\bibitem[9]{9} J. Chazy, 
{\em Sur les \'equations diff\'erentielles du trousi\'eme ordre et d'ordre sup\'erieur dont l'int\'egrale a ses points critiques fixes}, 
Acta Math. {\bf 34} (1911), 317--385. 



\bibitem[10]{Sasa1} Y. Sasano, {\em Coupled Painlev\'e VI systems in dimension four with affine Weyl group symmetry of types $B_6^{(1)},D_6^{(1)}$ and $D_7^{(2)}$}, preprint.

\bibitem[11]{Sasa2} Y. Sasano, {\em Four-dimensional Painlev\'e systems of types $D_5^{(1)}$ and $B_4^{(1)}$}, preprint.

\bibitem[12]{Sasa5} Y. Sasano, {\em Higher order Painlev\'e equations of type ${D_l}^{(1)}$}, RIMS Kokyuroku {\bf 1473} (2006), 143--163.

\bibitem[13]{Sasa4} Y. Sasano, {\em Symmetries in the system of type $D_4^{(1)}$}, preprint.

\bibitem[14]{Sasa3} Y. Sasano, {\em Coupled Painlev\'e III systems with affine Weyl group symmetry of types $B_4^{(1)}$, $D_4^{(1)}$ and $D_5^{(2)}$}, preprint.

\bibitem[15]{Sasa6} Y. Sasano, {\em Coupled Painlev\'e III systems with affine Weyl group symmetry of types $B_5^{(1)},D_5^{(1)}$ and $D_6^{(2)}$}, preprint.

\bibitem[16]{Sasa7} Y. Sasano, {\em Coupled Painlev\'e VI systems in dimension four with affine Weyl group symmetry of type $D_6^{(1)}$, II}, RIMS Kokyuroku Bessatsu. {\bf B5} (2008), 137--152.

\bibitem[17]{T} T. Tsuda, K Okamoto and H. Sakai, {\em Folding transformations of the Painlev\'e equations}, Math. Ann. {\bf 331} (2005), 713-738.

\end{thebibliography}
\end{document}